\documentclass[]{MathAppl18}
\renewcommand{\maketitle}{}
\usepackage[OT4,T1]{fontenc}
\usepackage[utf8]{inputenc}
\usepackage{graphicx}
\usepackage{csquotes}
\usepackage{amsmath}
\usepackage{amssymb}
\usepackage{filecontents}
\usepackage{lastpage}
\usepackage{enumerate}
\usepackage{lastpage} 

\usepackage{etex} 

\usepackage{color}
\usepackage{cite} 

\RequirePackage[numbers]{natbib}

\usepackage[colorlinks=true]{hyperref}
  \hypersetup{
    pdftitle={Teplate Mathematica Applicanda }, 
    pdfauthor=Nowy Autor, 
    colorlinks,
    urlcolor=blue,
    filecolor=magenta,
    citecolor=blue,
    linkbordercolor={1 1 1}, 
    citebordercolor={1 1 1},  
    urlbordercolor={ 1 1 1}  
  }

\usepackage{graphicx}
\usepackage{wrapfig}
\graphicspath{{./Figures/},{./Pictures/}}
\usepackage{multicol}

\numberwithin{equation}{section}
\firstpage{i}

\makeatletter
\renewcommand{\qed}{\@ifnextchar\par{\@gobble\qed}{\nobreak\hfill\hbox{\@square}\par}}
\makeatother

\title{\bf Complete Classification of the Euclidean Complete\\[10pt] Solutions to a Monge-Amp\`{e}re Equation}
\tpauthortrue
\author[Shi-Zhong Du et al.]{Shi-Zhong Du\footnotemark[1],~~Chen-Long Wu\footnotemark[1],~~Fei-Hao Zheng\footnotemark[1]}


\usepackage{marvosym}

\makeatletter
\def\@fnsymbol#1{%
  \ifcase#1\or
  *\or
  \Letter\fi
}
\makeatletter
\renewcommand{\@makefntext}[1]{%
  \hangindent=1.5em
  \noindent
  \makebox[1.5em][l]{\@thefnmark}%
  #1%
}
\makeatother

\begin{document}

\maketitle
\footnotetext[1]{{\scriptsize The Department of Mathematics, Shantou University, Shantou 515063, P.R. China}}

\renewcommand{\thefootnote}{\Letter}
\footnotetext[2]{{\scriptsize  S.Z. Du (szdu@stu.edu.cn), C.L. Wu (24clwu@stu.edu.cn), F.H. Zheng (20fhzheng@stu.edu.cn)}}\vspace{-5ex}
\renewcommand{\thefootnote}{\ensuremath{*}}
\footnotetext[3]{{\scriptsize The authors were partially supported by NFSC (12171299) and GDNSF (2019A1515010605).}}
\setcounter{footnote}{0}
\renewcommand{\thefootnote}{\arabic{footnote}}

\setcounter{page}{1}

\begin{abstract}
We study the Monge-Amp\`{e}re equation
     \begin{equation}\label{e0.1}
     \begin{aligned}
       \det D^2u=u^p, \ \ \forall x\in\Omega
     \end{aligned}
     \end{equation}
   for some $p\in{\mathbb{R}}$. A solution $u$ of \eqref{e0.1} is called to be Euclidean complete if it is an entire solution defined over the whole ${\mathbb{R}}^n$ or its graph is a large hypersurface satisfying the large condition $u(x)\to\infty$ as $\mathrm{dist}(x,\partial\Omega)\to 0$ in case of $\Omega\not={\mathbb{R}}^n$. In this paper, we will give various sharp conditions on $p$ and $\Omega$ classifying the Euclidean complete solution of \eqref{e0.1}.
\end{abstract}

\section{Introduction}

\bigbreak
\noindent The Dirichelt problem for the Monge-Amp\`{e}re equation
\begin{equation}\label{e1.1}
\det D^2u=f(u),~~\forall x\in\Omega\subset\mathbb{R}^n,
\end{equation}
and more generally for the $k$-Hessian equation $S_k(D^2u)=f(u)$ (with $\Delta u=S_1(D^2u)$ and $\det D^2u=S_n(D^2u)$ as special cases) has been studied extensively in recent decades. Among these works, several groundbreaking results are worth mentioning. For instance, Trudinger and Urbas \cite{TU} established the necessary and sufficient conditions for the solvability of \eqref{e1.1} with $f(u)$ replaced by the more general $f(x,u,Du)$.
In \cite{Wang1}, Wang studied the existence theory for classical solutions of Hessian equations and  derived the Hessian-Sobolev inequality
\begin{equation}\label{e1.2}                                                                                                                                             \left\|u\right\|_{L^p(\Omega)}^{k+1}\le -C\int_{\Omega}S_k(D^2u)udx,~~\forall p\in [1,k^*].
\end{equation}
Here $k^*$ denotes the critical exponent for the $k$-Hessian operator: $k^*=\frac{n(k+1)}{n-2k}$ when $2k<n$; $k^*$ can be any real number greater than $1$ when $2k=n$; and
$k^*=\infty$ whenever $2k>n$. Obviously, for $k=1$, inequality \eqref{e1.2} reduces to the classical Sobolev inequality. When $k=n$, the extremal functions of inequality \eqref{e1.2} satisfy \eqref{e1.1}. Conversely, the existence of solutions to \eqref{e1.1} can also be derived by means of the Hessian-Sobolev inequality and variational methods. In addition, Trudinger and Wang \cite{TW} provided a measure characterization of the
$k$-Hessian operator $S_k(D^2u)$. They proved that for any $k$-convex function $u\in \Phi^k(\Omega)$ (i.e., $u$ is upper semicontinuous in $\Omega$ and satisfies $S_k(D^2u)\ge 0$ in the viscosity sense), there exists a Borel measure $\mu_k[u]$ such that $\mu_k[u]=S_k(D^2u)$
for $u\in C^2(\Omega)$, and if $\{u_m\}$ is a sequence in $\Phi^k(\Omega)$ converging locally in measure to a function $u\in \Phi^k(\Omega)$, the sequence of Borel measures $\{\mu_k[u_m]\}$ converges weakly to $\mu_k[u]$. In \cite{CW1}, Chou and Wang systematically developed the variational theory of Hessian equations,
generalizing the well-known Ambrosetti-Rabinowitz theory for semilinear elliptic equations to the Hessian equations. For more results concerning the Dirichlet problems
for Monge-Amp\`{e}re and Hessian equations,
see \cite{Tru2, Tru3, CNS, Caf1, Guan, Lion1, Lion2, Pog, Tso, Tru1} and the references therein.

\medbreak
An interesting problem is whether \eqref{e1.1} admits large solutions or boundary blow-up solutions (i.e., satisfying $u=\infty$ on the boundary $\partial\Omega$) when the nonlinear terms $f(u)$ satisfies certain conditions. In fact, the study of boundary blow-up problems for elliptic partial differential equations can be traced back to the work of Bieberbach \cite{Bie1} and Rademacher \cite{Rad1} on the equation $\Delta u=e^u$. Let us review the following boundary blow-up problem involving the standard Laplacian
\begin{equation}\label{e1.3}
\Delta u=f(u),~~\forall x\in\Omega\subset\mathbb{R}^n,~~u|_{\partial\Omega}=\infty,
\end{equation}
where the boundary condition means that $u(x)\to\infty$ as $d(x):=\mathrm{dist}(x,\partial\Omega)\to 0$. When $f(u)=e^{u}$, and $\Omega$ is a smooth bounded domain in $\mathbb{R}^2$, Bieberbach proved that there is a solution $u\in C^2(\Omega)$ satisfying
\begin{equation*}
\left|u(x)-\log ((d(x))^{-2})\right|=o(1), ~~\mathrm{as}~~d(x)\to 0.
\end{equation*}
Later, Keller \cite{Kel1} and Osserman \cite{Oss1} first supplied the following necessary and sufficient condition for the existence of solutions to problem \eqref{e1.3}:
\begin{equation*}
\int_{u}^{\infty}\left(2F_0(s)\right)^{-\frac{1}{2}}ds <\infty,~~\forall u>0,
\end{equation*}
where $F_0(s):=\int_0^sf(s)ds$. If $f(u)=u^{\frac{n+2}{n-2}},~n>2$,  Loewner and Nirenberg \cite{LN1} showed that the problem \eqref{e1.3} has a unique positive solution $u$ which satisfies
\begin{equation*}
u(x)\left(d(x)\right)^{\frac{n-2}{2}}=\left(\frac{n(n-2)}{4}\right)^{\frac{n-2}{2}},~~\mathrm{as}~~d(x)\to 0.
\end{equation*}
Subsequently, their results were extended to $f(u)=u^p$ with $p>1$ by Kondrat\'{e}v and Nikishkin \cite{KN1}. Since then, the existence, uniqueness and asymptotic behavior of large solutions (or boundary blow-up solutions) for semilinear or quasilinear elliptic PDEs has been extensively studied, and see
\cite{BM1, BM2, CD, CN1, CR1,CR2, MV1, MV2, ZMML} and their references.

\medbreak
When $\Omega$ is a subdomain of $\mathbb{R}^n$, the large solutions to the weighted Monge-Amp\`{e}re equation
\begin{equation}\label{e1.4}
\det D^2u=b(x)f(u),~~x \in\Omega,~~u|_{\partial\Omega}=\infty
\end{equation}
has also attracted considerable attention from numerous scholars. Similary, the boundary condition in \eqref{e1.4} means that $u(x)\to\infty$ as $d(x)=\mathrm{dist}(x,\partial\Omega)\to 0$. The problem \eqref{e1.4} was first been considered and studied by Cheng and Yau \cite{CY1, CY2} for exponential nonlinear term $f(u)=e^{Ku}$ in bounded convex domain and $b(x)f(u)=e^{2u}$ in unbounded domain due to their applications in geometry.  They established that for $b(x)\in C^{k-2,\alpha}(\Omega)$, problem \eqref{e1.4} admits a strictly convex solution
$u\in C^{k,\alpha}(\Omega)$, which is unique if $b(x)$ is analytic in $\Omega$.
In \cite{LM1}, Lazer and McKenna showed that problem \eqref{e1.4} has a unique strictly convex solution $u\in C^{\infty}(\Omega)$ for $b\in C^{\infty}(\overline{\Omega})$ with $b(x)>0$ on $\Omega$, and $f(u)=e^u$ or $f(u)=u^p$ with $p>n$. Subsequently, Mohammed \cite{Moh1} provided two sufficient conditions for the existence and non-existence of large solutions on strictly convex bounded domain, respectively. In particular, his results were further improved by Zhang and Du in \cite{ZD}, who proved that under certain conditions, problem \eqref{e1.4} admits a strictly convex solution if and only if
\begin{equation*}
\int_t^{\infty}\left[(n+1)F(s)\right]^{-\frac{1}{n+1}}ds<\infty, ~~\forall t>\tau,
\end{equation*}
where $\tau\in [-\infty,\infty)$ is some extended real number, $F$ is the antiderivative of $f$ with $F(\tau)=0$.
For the special case where $f(u)=u^p$ or $f(u)=e^u$, Yang and Chang \cite{YC} investigated
in detail the asymptotic behavior of large solutions near the boundary. In \cite{Zhang4}, Zhang generalized the results of Yang and Chang and obtained the optimal
global and boundary asymptotic property of large solutions to problem \eqref{e1.4}. More related results can be found
in \cite{CT, Mat1, GJ, Zhang1, Pli, Zhang3, Sal1, CSF, Tak1, ZF, ML, Huang1}.
When the nonlinear term $f(u)$ in the Monge-Amp\`{e}re equation \eqref{e1.1} is replaced by the more general nonlinear gradient term $f(x,u,Du)$, the study of existence,
uniqueness, and asymptotic behavior near the boundary of large solutions becomes more complicated;  we refer the reader to \cite{FZ, Zhang2, Zhang5, ZF2}
and the references therein for further details.

\medbreak
The main purpose of this paper is to investigate the existence and non-existence of Euclidean complete solutions for the Monge-Amp\`{e}re equation \eqref{e1.1} in the case  that $\Omega$ is a bounded or unbounded domain in $\mathbb{R}^n$. Throughout this paper, a solution $u$ of \eqref{e1.1} is said to be Euclidean complete if it is either an
entire solution defined over the whole $\mathbb{R}^n$, or a large solution satisfying the boundary blow-up condition
\begin{equation}\label{e1.5}
\lim_{x\to \partial\Omega}u(x)=\infty
\end{equation}
on the boundary $\partial\Omega$ in the case where $\Omega\ne \mathbb{R}^n$.

\medbreak
When $\Omega\subset\mathbb{R}^n$ is bounded domain, a first satisfactory existence result was shown by Cr\^{i}stea-Trombetti in \cite{CT}. Using a hypothesis \[
\int_{1}^{\infty} \frac{du}{\sqrt{F(u)}} < \infty,
\quad \text{where } F(u) := \int_{0}^{u} f_0^{1/n}(s) \, ds
\leqno{(\mathcal{H}_1)}\]
together with another mild assumption
\[
\left\{
\begin{aligned}
&f_0^{1/n}~\text{is locally Lipschitz continuous on}~[0,\infty),\\
&\text{convex, positive and non-decreasing on}~(0,\infty)~\text{with}~f_0(0)=0,
\end{aligned}
\right.
\leqno{(\mathcal{H}_2)} \]
they have proven the following result. \begin{theorem}\label{t1.1}
Let $\Omega$ be a smooth, strictly convex, bounded domain in $\mathbb{R}^n$ with $n\ge 2$. Supposing that there exists a function $f_0$ on $[0,\infty)$ which satisfies
the hypotheses $(\mathcal{H}_1)$ and $(\mathcal{H}_2)$, if $f_0(u)\le f(u)$ for every $u>0$, then \eqref{e1.1} admits a strictly convex large solution on $\Omega$.
\end{theorem}

This result was later extended to Hessian equations by Huang in \cite{Huang1}, and further generalized to Monge-Amp\`{e}re equations with more general nonlinear terms  by Zhang, Du-Zhang in \cite{ZD, Zhang1, Zhang2, Zhang3, Zhang4}. The first part of this paper is devoted to proving the counterpart of Theorem \ref{t1.1} as follows.

\begin{theorem}\label{t1.2}
Supposing that $f$ is an unbounded function satisfying
\begin{equation}\label{e1.6}
\limsup_{u\to\infty}\frac{f(u)}{u^n}<\infty,
\end{equation}
then for any bounded smooth convex domain $\Omega\subset\mathbb{R}^n$, there is no convex large solution to \eqref{e1.1}.
\end{theorem}

Unlike the case of bounded domains, when considering the entire solution defined the whole of $\mathbb{R}^n$, the situation changes dramatically. In fact, when
$f(u)=u^p$ and $p>n$, the results of Jin, Li and Xu in \cite{JLX1} on the Hessian equation show that problem \eqref{e1.1} has no positive entire convex subsolution.
In \cite{JB1}, Ji and Bao showed that problem \eqref{e1.1} admits a positive radial entire large solution on $\mathbb{R}^n$ if and only if $0<p<n$. Recently, Li and Bao
\cite{LB1} obtained some new existence and nonexistence results for nonradial entire large solutions of the Hessian equation $S_k(D^2u)=b(x)u^p$ in the case $0<p<k$. For
other works related to entire solutions of the Monge-Amp\`{e}re or Hessian equations, we refer the reader to \cite{BCGJ, CW2, JW1, ZZ1} and the references therein. In
particular, for the case $p>n$, we present the following sharp non-existence result for entire solutions of \eqref{e1.1}.

\begin{theorem}\label{t1.3}
Supposing that $f\in C((0,\infty))$ is a positive monotone non-decreasing function on $u$ and satisfies
\begin{equation}\label{e1.7}
f(u)\ge \mathcal{A}u^p,~~\forall u>0
\end{equation}
for some $p>n$ and $\mathcal{A}>0$, then there is no positive entire convex solution of \eqref{e1.1} on $\mathbb{R}^n$.
\end{theorem}

Although Theorem \ref{t1.2} and Theorem \ref{t1.3} may have been known before, for the sake of completeness, we present two simple proofs here for the convenience
of the readers. With the help of the techniques developed in \cite{Du}, we can also prove the following sharp existence result for entire solutions in the case $p<n$.

\begin{theorem}\label{t1.4}
Considering $f(u)=u^p$ for $p<n$ and $\Omega=\mathbb{R}^n$ in \eqref{e1.1}, there are infinitely many positive convex solutions which are affine inequivalent. More
precisely, for any $a_0>0$, there exists at least one convex solution $u$ of \eqref{e1.1} satisfying
\begin{equation*}
u(0)=a_0,~~Du(0)=0.
\end{equation*}
\end{theorem}

Finally, we turn to the case of unbounded domains $\Omega\ne \mathbb{R}^n$ and establish the following non-existence result for $p>n$.

\begin{theorem}\label{t1.5}
Supposing that $f\in C((0,\infty))$ is positive non-decreasing function on $u$ and satisfies \eqref{e1.7} for some $p>n$, $n\ge 2$ and $\mathcal{A}>0$, there is no
positive convex solution of \eqref{e1.1} for ideal domain $\Omega$.
\end{theorem}

As illustrated in Theorems \ref{t1.3}-\ref{t1.5}, it is natural to conjecture that for $p<n$ and unbounded domain $\Omega\ne \mathbb{R}^n$, there exists some positive
convex solutions which are Euclidean complete. However, focusing on the case of $n=2$ and $p\in (0,1/2)$, we have shown the following surprising non-existence result.

\begin{theorem}\label{t1.6}
Considering $f(u)=u^p$ for $n=2, ~p\in (0,1/2)$ and unbounded domain $\Omega\ne \mathbb{R}^2$ in \eqref{e1.1}, there is no positive convex solutions which is Euclidean
complete.
\end{theorem}

The contents of this paper are organized as follows. We will prove Theorem \ref{t1.2} in Sect. 2, and prove Theorem \ref{t1.3} in Sect. 3. Using the techniques developed
in \cite{Du}, we prove Theorem \ref{t1.4} in Sect. 4-5. Finally, the proofs of Theorem \ref{t1.5} and Theorem \ref{t1.6} are presented in Sect. 6 and Sect. 7
respectively.
\bigbreak
\section{Non-existence of large solution on bounded domain}

\bigbreak
\noindent Let $\Omega$ be a bounded domain in $\mathbb{R}^n$, and consider the following Dirichlet problem
\begin{equation}\label{e2.1}
\begin{cases}
\det D^2u=f, & \mbox{ in } \Omega,\\
u=0, & \mbox{ on } \partial\Omega,
\end{cases}
\end{equation}
let us first recall an important lemma established by Caffarelli in \cite{Caf1}.

\begin{lemma}\label{l2.1}
Consider \eqref{e2.1} where $\Omega$ is a convex domain satisfying that $B_\sigma\subset\Omega\subset B_{\sigma^{-1}}$ for some positive constant $\sigma$ and $0\leq f(x)\leq C_1$ for any $x\in\Omega$. Then the convex solution $u$ of \eqref{e2.1} satisfies that
\begin{equation}\label{e2.2}
u(x)\geq -C_2\mathrm{dist}^\gamma(x,\partial\Omega),  \ \ \forall x\in\Omega,
\end{equation}
where $\gamma=\frac{2}{n}$ for $n\geq3$ and $\gamma\in(0,1)$ for $n=2$, $C_2=C_2(n,\sigma,C_1)$.
\end{lemma}

As a counterpart of Theorem \ref{t1.1}, we have the following non-existence result.

\begin{theorem}\label{t2.1}
Supposing that {$f$ is an unbounded function satisfying}
\begin{equation}\label{e2.3}
\limsup_{u\to\infty}\frac{f(u)}{u^n}<\infty,
\end{equation}
{then} for any bounded domain $\Omega\subset{\mathbb{R}}^n$, there is no convex large solution to \eqref{e1.1}.
\end{theorem}

\begin{proof}
By assumption \eqref{e2.3}, there exists an increasing sequence $u_k\to\infty$ such that
\begin{equation}\label{e2.4}
f(u_k):=\sup_{u\leq u_k}f(u), \ \ \lim_{k\to\infty}u_k^{-n}f(u_k)={\mathcal{A}}\in[0,\infty).
\end{equation}
Letting $\Sigma_{\mathcal{B}}:=\{x\in\Omega~|~\mathrm{dist}(x,\partial\Omega)>\mathcal{B}\}$, we claim that
\begin{equation}\label{e2.5}
u(x)=\infty, \ \ \forall x\in\Omega\setminus\Sigma_{{\mathcal{B}}}
\end{equation}
holds for any ${\mathcal{B}}<(C_2{\mathcal{A}}^{1/n})^{-1/\gamma}$. In fact, setting $\Omega_k:=\{x\in\Omega~|~u(x)<u_k\}$ and rescaling $v_k(x):=f^{-1/n}(u_k)u(x)$, for sufficiently large $k\in\mathbb{N}$, we have $\varnothing\ne \Omega_k\subset\Omega_{k+1}$ and
\begin{equation}\label{e2.6}
\det D^2v_k=f(u)/f(u_k)\leq 1, \ \ \forall x\in\Omega_k.
\end{equation}
Hence, by \eqref{e2.6} and Lemma \ref{l2.1}, it follows that
\begin{eqnarray*}
f^{-1/n}(u_k)u(x)-u_kf^{-1/n}(u_k)&\geq& -C_2\mathrm{dist}^\gamma(x,\partial\Omega_k)\\
&\geq&-C_2{\mathcal{B}}^\gamma
\end{eqnarray*}
for any $x\in\Omega_k\setminus\Sigma_{{\mathcal{B}}}$, this is equivalent to
\begin{equation}\label{e2.7}
u(x)\geq u_k-C_2{\mathcal{B}}^\gamma f^{1/n}(u_k), \ \ \forall x\in\Omega_k\setminus\Sigma_{{\mathcal{B}}}.
\end{equation}
Sending $k\to\infty$ in \eqref{e2.7} and using \eqref{e2.4} for ${\mathcal{B}}<(C_2{\mathcal{A}}^{1/n})^{-1/\gamma}$, we conclude that $u$ is infinite everywhere in $\Omega\setminus\Sigma_{{\mathcal{B}}}$. The claim \eqref{e2.5} holds true. By an bootstrapping argument we reach the conclusion that $u$ is infinite everywhere in the whole domain $\Omega$. This completes the proof of the theorem.
\end{proof}

\bigbreak \section{Non-existence of entire solution on ${\mathbb{R}}^n$}

\bigbreak	
\noindent In this section, we will show the following sharp non-existence result for entire solutions to the Monge-Amp\`{e}re equation \eqref{e1.1}.
\begin{theorem}\label{t3.1}
Supposing that $f\in C((0,\infty))$ is a positive monotone non-decreasing function on $u$ and satisfies
\begin{equation}\label{e3.1}
f(u)\geq{\mathcal{A}}u^p, \ \ \forall u>0
\end{equation}
for some $p>n$ and ${\mathcal{A}}>0$, then there is no positive entire convex solution of \eqref{e1.1} in ${\mathbb{R}}^n$.
\end{theorem}

Before proving the theorem, we need the following existence result by Lazer-McKenna in \cite{LM1}.

\begin{lemma}\label{l3.1}
Letting $\Omega$ be a bounded smooth convex domain in ${\mathbb{R}}^n$ and $f(u)=u^p$ for some $p>n$, there exists an unique large convex solution $u\in C^\infty(\Omega)$ of \eqref{e1.1} satisfying
\begin{equation}\label{e3.2}
C_3^{-1}\mathrm{dist}^{-\alpha}(x,\partial\Omega)\leq u(x)\leq C_3\mathrm{dist}^{-\alpha}(x,\partial\Omega)
\end{equation}
for some positive constant $C_3=C_3(n,p,\Omega)$, where $\alpha=\frac{n+1}{p-n}$.
\end{lemma}

\noindent\textbf{Proof of Theorem \ref{t3.1}:}  Let $u$ be a positive entire convex solution of \eqref{e1.1} in $\mathbb{R}^n$. Then, by assumption \eqref{e3.1}, $u$ is also a subsolution of
\begin{equation}\label{e3.3}
\det D^2u={\mathcal{A}}u^p, \ \ \forall x\in\Omega,
\end{equation}
where $\Omega$ is an arbitrary open domain in $\mathbb{R}^n$. By Lemma \ref{l3.1}, we denote by $u_{\lambda}$ the unique convex large solution of \eqref{e1.1} with $f(u)=u^p$ and $\Omega=B_{\lambda}$, and rescale $u_{\lambda}$ by
$$
v_{\lambda}(x):=\lambda^{\frac{2n}{p-n}}u_{\lambda}(\lambda x), \quad \forall x\in B_1.
$$
It is clear that $v_{\lambda}$ is the unique convex large solution of
$$
\begin{cases}
\det D^2v_{\lambda}=v_{\lambda}^p, ~~&\text{in}~~B_1,\\
v_{\lambda}(x)\to \infty,~~&\text{as}~~x\to \partial B_1.
\end{cases}
$$
Applying Lemma \ref{l3.1} again, we obtain that
$$
C_3^{-1}\mathrm{dist}^{-\alpha}(x,\partial B_1)\leq v_{\lambda}(x)\leq C_3\mathrm{dist}^{-\alpha}(x,\partial B_1), \quad \forall x\in B_1.
$$
Rewriting this estimate in terms of the function $u_{\lambda}(y)$, we have
\begin{equation}\label{e3.4}
C_3^{-1}\lambda^{\frac{2n}{n-p}+\alpha}\mathrm{dist}^{-\alpha}(y,\partial B_{\lambda})\leq u_{\lambda}(y)\leq
C_3\lambda^{\frac{2n}{n-p}+\alpha}\mathrm{dist}^{-\alpha}(y,\partial B_{\lambda})
\end{equation}
which holds for arbitrary $y\in B_{\lambda}$. Defining $w_{\lambda}(y):={\mathcal{A}}^{\frac{1}{n-p}}u_{\lambda}(y)$, then $w_{\lambda}$ is a solution to
$$
\begin{cases}
\det D^2w_{\lambda}={\mathcal{A}}w_{\lambda}^p, ~~&\text{in}~~B_{\lambda},\\
w_{\lambda}(y)\to\infty,~~&\text{as}~~y\to\partial B_{\lambda}.
\end{cases}
$$
Next, let us denote by $\left[U^{ij}(t)\right]$ the cofactor matrix of $tD^2w_{\lambda}+(1-t)D^2u$ with $t\in [0,1]$. By some simple calculations, it is no difficult to find that
$$
\begin{cases}
\displaystyle \int^1_0U^{ij}(t)dt\cdot D_{ij}(w_{\lambda}-u)\le {\mathcal{B}}(w_{\lambda}-u), & \text{in}~~B_{\lambda},\\
w_{\lambda}-u\geq0, & \text{on}~~\partial B_{\lambda},
\end{cases}
$$
where the coefficient $\mathcal{B}$ is defined as
$$
{\mathcal{B}}=\frac{f(w_{\lambda})-f(u)}{w_{\lambda}-u}\geq 0.
$$
It is inferred from the maximum principle of elliptic PDEs that $u\leq w_{\lambda}$ in $B_{\lambda}$. Combining \eqref{e3.4} and the definition of $w_{\lambda}$, we immediately obtain
\begin{equation}\label{e3.5}
u(y)\leq C_3{\mathcal{A}}^{\frac{1}{n-p}} \lambda^{\frac{2n}{n-p}+\alpha}\mathrm{dist}^{-\alpha}(y,\partial B_{\lambda}),\quad \forall y\in B_{\lambda}.
\end{equation}
Fixing $y$ and letting $\lambda$ tends to infinity, it is inferred from \eqref{e3.5} that $u$ is identical to zero for $p>n$. This contradicts the assumption that $u$ is the positive solution, and thus the proof is complete.\hfill\ensuremath{\square}

\bigbreak

\section{Schauder's fixed point scheme and entire solutions for $p<n$}

\bigbreak
\noindent Letting $u=u(r)$ be a radial symmetric solution to \eqref{e1.1} with $f(u)=u^p$, there holds that
\begin{equation}\label{e4.1}
u''(u')^{n-1}r^{1-n}=u^p.
\end{equation}
Noting that for $p<n$,
$$
u(x)=\beta_{n,p}|x|^\alpha, \ \ \alpha:=\frac{2n}{n-p},\ \beta_{n,p}:=\Big[\alpha^n(\alpha-1)\Big]^{\frac{1}{p-n}}
$$
is a nonnegative convex entire solution of \eqref{e1.1} for $f(u)=u^p$ with one degenerate singular point $x=0$, it can not be expect to prove an A-priori positive lower bound and regularity for entire solution of \eqref{e1.1} as that in \cite{CT} for large solution. However, one may also ask the following question of existence of positive strictly convex entire solution of \eqref{e1.1}.

\medbreak
\noindent\textbf{Question.} Considering \eqref{e1.1} for $f(u)=u^p, ~p<n$ and $\Omega={\mathbb{R}}^n$, whether it admits a smooth positive convex entire solution?

\medbreak
We will give a confirm answer in the following theorem.

\begin{theorem}\label{t4.1}
Considering \eqref{e1.1} for $f(u)=u^p, ~p<n$ and $\Omega={\mathbb{R}}^n$, there are infinitely many positive entire convex solutions which are affine inequivalent. More precisely, for any $a_0>0$, there exists at least one positive entire convex solution $u$ of \eqref{e1.1} satisfying
$$
u(0)=a_0, \ \ Du(0)=0.
$$
\end{theorem}

Before proving the theorem, let us check an easy necessary condition for positive entire solution.

\begin{lemma}\label{l4.1}
A necessary condition ensuring
\begin{equation}\label{e4.2}
u(r)=\sum_{j=0}^{\infty}\frac{a_j}{j!}r^j, \ \ \forall r\in[0,\delta]
\end{equation}
to be a radial symmetric positive entire solution is given by
\begin{equation}\label{e4.3}
a_0, a_2>0, \ \ a_2=a_0^{p/n}, \ \ a_{2j-1}=0, \ \ \forall j\in{\mathbb{N}}.
\end{equation}
\end{lemma}

\begin{proof}
The lemma can be verified by substituting $u$ into \eqref{e4.1} at $r=0$.
\end{proof}

At the first step, we introduce an iteration scheme to prove a local existence result for \eqref{e4.1}. Letting $\kappa\in\mathbb{N}$, $a:=(a_0,a_1,\cdots,a_{2\kappa})\in{\mathbb{R}}^{2\kappa+1}$ be a vector satisfying \eqref{e4.3} and $\sigma,\delta\in(0,1)$, define $\mathcal{F}_{a,\delta,\sigma,\kappa}$ to be the set consisting of all functions $u\in C^{2\kappa}([0,\delta])$ such that
\begin{equation*}
\left\{
\begin{aligned}
&u^{(j)}(0)=a_{j},~u^{(j)}(r)\in[a_{j}-\sigma,a_{j}+\sigma],\\
&u^{(2\kappa)}(0)=a_{2\kappa},~u^{(2\kappa)}(r)\in[a_{2\kappa}-1,a_{2\kappa}+1]
\end{aligned}
\right.
\end{equation*}
for any $r\in [0,\delta],~j=1,2,\dots,2\kappa-1$. In what follows, $\mathcal{F}_{a,\delta,\sigma,\kappa}$ is considered a convex subset of $C^{2\kappa-1}([0,\delta])$ endowed with norm $\|\cdot\|_{C^{2\kappa-1}([0,\delta])}$. For any function $\varphi\in\mathcal{F}_{a,\delta,\sigma,\kappa}$, we define $\xi:=T\varphi$ to be the solution determined by
\begin{equation}\label{e4.4}
\left\{
\begin{aligned}
&\displaystyle\xi''(r)=\varphi^p(\varphi')^{1-n}r^{n-1},~r\in[0,\delta],\\
&\xi(0)=a_0,~\xi'(0)=0,~\xi''(0)=a_2.
\end{aligned}
\right.
\end{equation}

\medbreak
Next, we will show that $T$ is a continuous and compact mapping from $\mathcal{F}_{a,\delta,\sigma,\kappa}$ to $\mathcal{F}_{a,\delta,\sigma,\kappa}$. Thus, it is inferred from the following variant of Schauder's fixed point theorem that $T$ admits a fixed point in the closure of $\mathcal{F}_{a,\delta,\sigma,\kappa}$.

\begin{lemma}\label{l4.2} {\rm(Schauder's fixed point lemma)}
Letting $C$ be a convex subset of Banach space $X$ and $T:X\to X$ be a continuous and compact mapping from $C$ to $C$, then there exists a fixed point of $T$ on the
closure of $C$.
\end{lemma}

\begin{proof}
The proof of Lemma \ref{l4.2} is straightforward. In fact, by the continuity of $T$, it is easy to see that $T(\overline{C})\subset \overline{T(C)}\subset \overline{C}$.
Let $U\subset\overline{C}$ be an arbitrary bounded set. For each $y\in U$, choose a sequence $\{y_n\}\subset C$ such that $y_n\to y$.
Letting $U_0:=\left\{y_n~|~y_n\to y,~y\in U\right\}$, then $U_0\subset C$ is also bounded. By compactness of $T:C\to C$,
$\overline{T(U_0)}$ is compact in $X$. Since $T(y)\in \overline{T(U_0)}$ for every $y\in U$, we have $T(U)\subset\overline{T(U_0)}$.
The closure $\overline{T(U)}$ is a closed subset of the compact $\overline{T(U_0)}$, hence compact.
Thus, $T$ is also a continuous compact mapping from $\overline{C}$ into $\overline{C}$. In view of the classical Schauder's fixed point theorem, there exists a fixed point of $T$ in $\overline{C}$.
\end{proof}

The proofs of Theorem \ref{t4.1} were divided into several lemmas. For any positive function $f\in C^k([0,\delta])$, we have $f^{-1}\in C^k([0,\delta])$. Therefore, supposing that
$$
f^{(l)}(0)=:f_l, \quad l=0,1,\cdots,k,
$$
let us introduce a new notation for conjugate indices
$$
\overline{f}_l:=\overline{f}_l(f_0,\cdots,f_l):=\frac{d^lf^{-1}}{dr^l}\Big|_{r=0}\in{\mathbb{R}}
$$
which is determined by $f_i, i=0,1,\cdots,l$ for each $l\leq k$. We need the following lemma.

\begin{lemma}\label{l4.3}
Letting $\varphi\in C^k([0,\delta])$ be a positive function satisfying
\begin{equation}\label{e4.5}
\varphi(0)=\varphi_0=0,~~\frac{d^l\varphi}{dr^l}\Big|_{r=0}=\varphi_l\in{\mathbb{R}}, ~~l=1,2,\cdots,k,
\end{equation}
then function $\phi\in C^{k}([0,\delta])$ defined by
$$
\phi(r):=\begin{cases}
\displaystyle  \frac{r}{\varphi(r)}, & r\in(0,\delta]\\
\displaystyle  \frac{1}{\varphi_1}, & r=0
\end{cases}
$$
satisfies that
\begin{equation}\label{e4.6}
\overline{\phi}_l:=\frac{d^l\varphi^{-1}}{dr^l}\Big|_{r=0}=\frac{\varphi_{l+1}}{l+1},~~l=0,1,\cdots, k-1.
\end{equation}
\end{lemma}

\begin{proof}
Letting $f:=\frac{\varphi}{r}$, by Taylor's expansion, we have
\begin{equation*}
\begin{aligned}
\varphi(r)&=\sum_{j=1}^{k-1}\frac{\varphi_j}{j!}r^j+\int^r_0\frac{\varphi^{(k)}(t)}{(k-1)!}(r-t)^{k-1}dt,\\
f(r)&=\sum_{j=1}^{k-1}\frac{\varphi_j}{j!}r^{j-1}+\int^r_0\frac{\varphi^{(k)}(t)}{r(k-1)!}(r-t)^{k-1}dt.
\end{aligned}
\end{equation*}
Repeated differentiation of $f$ gives
\begin{equation}\label{e4.7}
\begin{aligned}
\frac{d^lf}{dr^l}=&\sum_{j=l+1}^{k-1}\frac{\varphi_j}{j(j-l-1)!}r^{j-l-1}+
\\&\sum_{q=0}^l\frac{(-1)^qC_l^qq!}{(k-l+q-1)!r^{q+1}}\int^r_0\varphi^{(k)}(t)(r-t)^{k-l+q-1}dt
\end{aligned}
\end{equation}
for $l=1,2,\cdots,k-2$ and
\begin{equation}\label{e4.8}
\frac{d^{k-1}f}{dr^{k-1}}=\sum_{q=0}^{k-1}\frac{(-1)^qC_{k-1}^q}{r^{q+1}}\int^r_0\varphi^{(k)}(t)(r-t)^qdt.
\end{equation}
By letting $r\to 0^+$ in \eqref{e4.7} and using Cauchy's Theorem, it follows that $f\in C^{k-2}([0,\delta])$ and
\begin{equation}\label{e4.9}
f_l:=\frac{d^lf}{dr^l}\Big|_{r=0}=\frac{\varphi_{l+1}}{l+1}, \quad l=0,1,\cdots,k-2.
\end{equation}
To show $f\in C^{k-1}([0,\delta])$, one needs only to use \eqref{e4.8} to deduce that
$$
\lim_{r\to0^+}\frac{d^{k-1}f}{dr^{k-1}}(r)=\lim_{\xi\to 0^+}\varphi^{(k)}(\xi)\sum_{q=0}^{k-1}\frac{(-1)^qC_{k-1}^q}{q+1}=\frac{\varphi_k}{k}
$$
by the integral mean value theorem, and thus conclude that $f\in C^{k-1}([0,\delta])$ and
\begin{equation}\label{e4.10}
f_{k-1}:=\frac{d^{k-1}f}{dr^{k-1}}\Big|_{r=0}=\frac{\varphi_k}{k}.
\end{equation}
Thus, \eqref{e4.6} follows from combining of \eqref{e4.9} and \eqref{e4.10}.
\end{proof}

\begin{lemma}\label{l4.4}
For any positive function $f\in C^k([0,\delta])$, one has
\begin{equation}\label{e4.11}
\overline{f}_l:=\frac{d^lf^{-1}}{dr^l}\Big|_{r=0}=-\frac{f_l}{f_0^2}+\Psi_l(f_0,f_1,\cdots,f_{l-1}),\quad l=1,\cdots,k
\end{equation}
for some rational functions $\Psi_l$ on $(f_0,f_1,\cdots, f_l)$.
\end{lemma}

The conclusion of the lemma follows from an easy calculation and induction on $l$. As a corollary of Lemma \ref{l4.3} and \ref{l4.4}, one gets that

\begin{corollary}\label{c4.1}
Under the assumptions of Lemma \ref{l4.3}, there holds
\begin{equation}\label{e4.12}
\phi_0=\frac{1}{\varphi_1},~~\phi_l=-\frac{\varphi_{l+1}}{(l+1)\varphi_1^2}+\Psi_l(\varphi_0,\varphi_1,\cdots,\varphi_l)
\end{equation}
for $l=1,2,\cdots,k-1$ and some rational function $\Psi_l$ on $(\varphi_0,\varphi_1,\cdots, \varphi_l)$.
\end{corollary}

Given $\varphi\in\mathcal{F}_{a,\delta,\sigma,\kappa}$ for some $a\in{\mathbb{R}}^{2\kappa+1}$ satisfying \eqref{e4.3}, the function $\xi=T\varphi$ satisfies that
\begin{equation}\label{e4.13}
\xi''(r)=\varphi^p(\varphi')^{1-n}r^{n-1}=\varphi^p\phi^{n-1},
\end{equation}
where
\begin{equation*}
\phi(r):=\left\{
\begin{aligned}
&\frac{r}{\varphi'(r)}, ~~&&r\in(0,\delta],\\
&\frac{1}{\varphi_2}, ~~&&r=0.
\end{aligned}
\right.
\end{equation*}
Using Corollary \ref{c4.1} for $\varphi$ replaced by $\varphi'$ and taking the $l$-th derivative of \eqref{e4.13} with respect to $r$ for $l=1,2,\cdots,2\kappa-2$, it yields that
\begin{equation*}
\begin{aligned}
\frac{d^{l+2}\xi}{dr^{l+2}}\Big|_{r=0}&=\sum_{q=0}^lC_l^q\frac{d^q\varphi^p}{dr^q}\frac{d^{l-q}\phi^{n-1}}{dr^{l-q}}\Big|_{r=0}\\
&=(n-1)\varphi_0^p\phi_0^{n-2}\left(-\frac{\varphi_{l+2}}{(l+1)\varphi_2^2}\right)+\widetilde{\Psi}_{l}(a_0,a_1,\cdots,a_{l+1})\\
&=-\frac{(n-1)a_0^pa_2^{-n}}{l+1}\frac{d^{l+2}\varphi}{dr^{l+2}}\Big|_{r=0}+\widetilde{\Psi}_{l}(a_0,a_1,\cdots,a_{l+1})
\end{aligned}
\end{equation*}
for some rational function $\widetilde{\Psi}_{l}$ on $(\varphi_0,\varphi_1,\cdots, \varphi_{l+1})$. It is remarkable that for $a\in{\mathbb{R}}^{2\kappa+1}$ satisfying \eqref{e4.3}, there holds
\begin{equation}\label{e4.14}
\frac{d^{2l-1}\xi}{dr^{2l-1}}\Big|_{r=0}=0,\quad l=1,2,\cdots,\kappa.
\end{equation}
Moreover, by taking higher derivatives on \eqref{e4.13} and applying the Cauchy's mean value theorem repeatedly, it is not hard to see that the following proposition holds true.

\begin{proposition}\label{p4.1}
Given $\varphi\in\mathcal{F}_{a,\delta,\sigma,\kappa}$ for some $a\in{\mathbb{R}}^{2\kappa+1}$ satisfying \eqref{e4.3}, the function $\xi=T\varphi$ satisfies that
\begin{equation}\label{e4.15}
\frac{d^{l+2}\xi}{dr^{l+2}}\Big|_{r=0}=-\mathcal{A}_l\frac{d^{l+2}\varphi}{dr^{l+2}}\Big|_{r=0}+\widetilde{\Psi}_{l}(a_0,a_1,\cdots,a_{l+1})
\end{equation}
and
\begin{equation}\label{e4.16}
\begin{aligned}
\frac{d^{l+2}\xi}{dr^{l+2}}=&-\left(\mathcal{A}_{l}+\Phi_{l,\varphi}\right)\frac{d^{l+2}\varphi}{dr^{l+2}}+\widetilde{\Psi}_{l}(a_0,a_1,\cdots,a_{l+1})+\widetilde{\Phi}_{l,\varphi}
\end{aligned}
\end{equation}
for $r\in[0,\delta]$ and $l=1,2,\cdots,2\kappa-2$, where $\mathcal{A}_l:=\left[(n-1)a_0^pa_2^{-n}\right]/(l+1)$ and the functions
$$
\left|\Phi_{l,\varphi}\right|\leq \sigma_1(l,a,\sigma),\quad \big|\widetilde{\Phi}_{l,\varphi}\big|\leq \sigma_2(l,a,\sigma)
$$
are small as long as $\sigma$ is small for each fixed $a\in{\mathbb{R}}^{2\kappa+1}$.
\end{proposition}

\noindent\textbf{Remark.} Without special indication, the functions on left hand side of \eqref{e4.16} take values on $r$, while the functions on right hand side of take values on $0<\cdots<r_3<r_2<r_1<r$, which come from the Cauchy's mean value theorem.

\medbreak
Using \eqref{e4.15} successively, then for any given $a_0, a_2$ satisfying \eqref{e4.3}, one can uniquely determine the whole vector $a$ by iterating
\begin{equation}\label{e4.17}
a_{l+2}=-\frac{(n-1)a_0^pa_2^{-n}}{l+1}a_{l+2}+\widetilde{\Psi}_{l}(a_0,a_1,\cdots,a_{l+1}),
\end{equation}
where $l=1,2,\cdots,2\kappa-2$. As mentioned above by \eqref{e4.14}, the whole vector $a$ satisfies also \eqref{e4.3}. Thus, we arrive at the following proposition.

\begin{proposition}\label{p4.2}
For some vector $a\in{\mathbb{R}}^{2\kappa+1}$ determined by \eqref{e4.17} and $a_0, a_2$, if one chooses $\kappa\in{\mathbb{N}}$ so large that
\begin{equation}\label{e4.18}
\mathcal{A}_{2\kappa-2}=\frac{(n-1)a_0^pa_2^{-n}}{2\kappa-1}<1
\end{equation}
and then $\sigma$ small, finally $\delta$ small, the mapping $T$ is a continuous and compact mapping from $\mathcal{F}_{a,\delta,\sigma,\kappa}$ to $\mathcal{F}_{a,\delta,\sigma,\kappa}$.
\end{proposition}

\begin{proof}
By virtue of \eqref{e4.16}, we immediately obtain that
\begin{equation}\label{e4.19}
\begin{aligned}
\xi^{(2\kappa)}=&-\left(\mathcal{A}_{2\kappa-2}+\Phi_{2\kappa-2,\varphi}\right)\varphi^{(2\kappa)}\\
&+\widetilde{\Psi}_{2\kappa-2}(a_0,a_1,\cdots,a_{2\kappa-1})+\widetilde{\Phi}_{2\kappa-2,\varphi}.
\end{aligned}
\end{equation}
Another hand, by iterative formula \eqref{e4.17} of $a$, one also has
\begin{equation}\label{e4.20}
a_{2\kappa}=-\mathcal{A}_{2\kappa-2}a_{2\kappa}+\widetilde{\Psi}_{2\kappa-2}(a_0,a_1,\cdots,a_{2\kappa-1}).
\end{equation}
Subtracting \eqref{e4.20} from \eqref{e4.19}, it follows that
\begin{equation}\label{e4.21}
\begin{aligned}
\left(\xi^{(2\kappa)}-a_{2\kappa}\right)=&-\left(\mathcal{A}_{2\kappa-2}+\Phi_{2\kappa-2,\varphi}\right)\left(\varphi^{(2\kappa)}-a_{2\kappa}\right)\\
&+\widetilde{\Phi}_{2\kappa-2,\varphi}-a_{2\kappa}\Phi_{2\kappa-2,\varphi}.
\end{aligned}
\end{equation}
From \eqref{e4.18}, if we choose $\sigma$ sufficiently small such that
\begin{equation}\label{e4.22}
\begin{aligned}
\vartheta:=\mathcal{A}_{2\kappa-2}&+\sigma_1(2\kappa-2,a,\sigma)\in(0,1),\\
|\sigma_2(2\kappa-2,a,\sigma)|& +|a_{2\kappa}\sigma_1(2\kappa-2,a,\sigma)|<1-\vartheta,
\end{aligned}
\end{equation}
it is inferred from \eqref{e4.21} that
\begin{equation}\label{e4.23}
\xi^{(2\kappa)}(r)\in[a_{2\kappa}-1,a_{2\kappa}+1]
\end{equation}
Therefore, if $\delta$ is chosen small, one has
\begin{equation}\label{e4.24}
\begin{aligned}
\xi^{(j)}(0)&=a_{j}, \ \ \xi^{(j)}(r)\in[a_{j}-\sigma,a_{j}+\sigma],\\
j&=0,1,\cdots,2\kappa-1, \ r\in[0,\delta].
\end{aligned}
\end{equation}
Combining \eqref{e4.23} and \eqref{e4.24}, we conclude that $T$ is a continuous and compact mapping from $\mathcal{F}_{a,\delta,\sigma,\kappa}$ to $\mathcal{F}_{a,\delta,\sigma,\kappa}$. The proof is complete.
\end{proof}

As a corollary, we have the following local existence result of positive entire solution of \eqref{e4.1}.

\begin{corollary}\label{c4.2}
Under assumptions of Proposition \ref{p4.2}, there exists a local smooth positive entire solution $u\in \mathcal{F}_{a,\delta,\sigma,\kappa}$ of \eqref{e4.1} on $[0,\delta]$ satisfying
\begin{equation}\label{e4.25}
u^{(j)}(0)=a_j, \quad j=0,1,\cdots,2\kappa
\end{equation}
for each $a$ determined by \eqref{e4.17} and $a_2=a_0^{p/n}$.
\end{corollary}

\bigbreak
\section{Long time existence and validity of Theorem \ref{t4.1}}

\bigbreak
\noindent Letting $u$ be the solution derived in Corollary \ref{c4.2}, we will show that the solution exists for all $r>0$, provided $p<n$.

\begin{lemma}\label{l5.1}
Letting $u$ be a solution of \eqref{e4.1} derived in Corollary \ref{c4.2} on $[0,R]$ for some $R>0$, if there holds
\begin{equation}\label{e5.1}
\sup_{r\in[0,R]}u^p(r)\leq C_4<\infty
\end{equation}
for some positive constant $C_4$, we have
\begin{equation}\label{e5.2}
\sup_{r\in[0,R]}\left[u'(r)+u''(r)\right]\leq C_5<\infty
\end{equation}
for positive constant $C_5=C_5(C_4,n,p,R,a_0,\delta,\sigma,\kappa)$.
\end{lemma}

\begin{proof}
Noting that by \eqref{e4.1}, the solution $u$ satisfies $u''>0$ and $u'>0$ on the maximal interval of existence. Since $u'$ is monotonically increasing and $u\in\mathcal{F}_{a,\delta,\sigma,\kappa}$, one has
\begin{equation}\label{e5.3}
\begin{aligned}
u'(r)&\geq u'(\delta)=\int^\delta_0u''(s)ds\\
&\geq(a_2-\sigma)\delta>0
\end{aligned}
\end{equation}
for any $r\in [\delta, R]$. Hence, it follows from \eqref{e4.1} and \eqref{e5.3} that
\begin{equation*}
\begin{aligned}
u''(r)&=u^p(u')^{1-n}r^{n-1}\\
&\leq C_4R^{n-1}[(a_2-\sigma)\delta]^{1-n}\\
&=:C_6<\infty
\end{aligned}
\end{equation*}
for arbitrary $r\in [\delta,R]$. In addition, by definition of $\mathcal{F}_{a,\delta,\sigma,\kappa}$, we have
$$
u''(r)\leq a_2+\sigma<\infty, \quad \forall r\in[0,\delta].
$$
That is to say, $u''$ is bounded from above by an A-priori bound $C_7$ depending only on $C_4,n,p,R,$ $a_0,\delta,\sigma,\kappa$. Integrating $u''$ over $[0,r]$, we have
$$
u'(r)=\int^r_0u''(s)ds\leq C_7R,~~\quad \forall r\in[0,R].
$$
This completes the proof of \eqref{e5.2}. \end{proof}

\begin{lemma}\label{l5.2}
Let $u$ be a solution of \eqref{e4.1} derived in Corollary \ref{c4.2}. If $p<n$, there exists a positive constant $C_8=C_8(n,p,R,a_0,\delta,\sigma,\kappa)$ such that
\begin{equation}\label{e5.4}
u^p(r)\leq C_8, \quad \forall r\in[0,R]
\end{equation}
holds for any $R>0$.
\end{lemma}

\begin{proof}
The case $p\leq 0$ is easy since the lower bound of $u$ follows from the monotonicity of $u$ and definition of $\mathcal{F}_{a,\delta,\sigma,\kappa}$. To solve the case $p\in(0,n)$, let us rewrite \eqref{e4.1} by
\begin{equation*}
\begin{aligned}
u''&=u^p(u')^{1-n}r^{n-1}\\
&\leq R^{n-1}u^p(u')^{1-n}.
\end{aligned}
\end{equation*}
Multiplying both sides of the above inequality by $(u')^n$, it yields that
\begin{equation*}
\frac{d}{dr}\left[\frac{1}{n+1}(u')^{n+1}-\frac{R^{n-1}}{p+1}u^{p+1}\right]\leq0.
\end{equation*}
After integrating over $[\delta,R]$, we conclude that
$$
u'\leq C_9(u+1)^{\frac{p+1}{n+1}}, \quad \forall r\in[\delta,R]
$$
for some positive constant $C_9=C_9(n,p,R,a_0,\delta,\sigma,\kappa)$. Thus, one gets
$$
\left[u(r)+1\right]^{\frac{n-p}{n+1}}\leq \frac{C_9(n-p)}{n+1}(R-\delta)+C_{10}
$$
for any $r\in [\delta,R]$ and some positive constant $C_{10}=C_{10}(n,p,R,a_0,\delta,\sigma,\kappa)$. Estimate \eqref{e5.4} holds since $p\in(0,n)$ and $u$ is A-priori bounded on $[0,\delta]$.
\end{proof}

\noindent\textbf{Complete the proof of Theorem \ref{t4.1}.} Lemma \ref{l5.1} and \ref{l5.2} guarantee that local solutions $u$ of \eqref{e4.1} obtained from Corollary \ref{c4.2} do not blow-up in finite time $r$, and hence yield the entire positive solutions of \eqref{e1.1}. Moreover, these solutions are affine inequivalent.\hfill\ensuremath{\square}

\bigbreak
\section{Large solutions on unbounded domain $\Omega\not={\mathbb{R}}^n$}

\bigbreak
\noindent In this section, we assume that $\Omega\not={\mathbb{R}}^n$ is an unbounded convex domain, which satisfies the following hypothesis of \enquote{ideal domain}:

\medbreak
\noindent $(\mathcal{H}_3)$ There exist a non-empty bounded Lipschitz convex domain $\Omega_0\subset{\mathbb{R}}^n$ and sequences $a_k\in\Omega, \lambda_k\to\infty$ such that
\begin{equation}\label{e6.1}
\Omega_k:=\lambda_k(\Omega_0-a_k)\subset\Omega, \ \ \forall k\in{\mathbb{N}}
\end{equation}
and
\begin{equation}\label{e6.2}
\inf_{k\in{\mathbb{N}}}\mathrm{dist}(a_\infty,\partial\Omega_k)>0
\end{equation}
holds for some $a_\infty\in\bigcap\limits_{k=1}^\infty\Omega_k$. We will call $\Omega$ to be an \enquote{ideal domain} and call $a_\infty$ to be an \enquote{ideal center} of $\Omega$.

\medbreak
\noindent\textbf{Remark.} Regarding the ideal domain, we provide the following remarks.
\medbreak
\noindent $(\mathcal{R}_1)$ If the bounded convex Lipschitz domain $\Omega_0$ in the definition of ideal domain is replaced by a ball $B$, then many unbounded convex domains will no longer satisfy the definition of ideal domain; a typical example is the first quadrant in $\mathbb{R}^2$.

\medbreak
\noindent $(\mathcal{R}_2)$ There exist unbounded convex domains that are not ideal domains, such as the ‘‘strip domain'' bounded by two parallel straight lines on the plane.

\medbreak
\noindent $(\mathcal{R}_3)$ All unbounded convex sets $\Omega\subset{\mathbb{R}}^n$ containing an infinite cone
$$
{\mathcal{C}}:=\Big\{\lambda x\in{\mathbb{R}}^n~\Big|\ x\in{\mathcal{C}}_0, \lambda\in{\mathbb{R}}^+\Big\},
$$
which is defined for a nonempty domain ${\mathcal{C}}_0\subset{\mathbb{R}}^n$, are ideal domains.

\medbreak
Next, we will prove the following nonexistence result for ideal domain.

\begin{theorem}\label{t6.1}
Supposing that $f\in C((0,\infty))$ is a positive monotone non-decreasing function on $u$ and satisfies \eqref{e3.1} for some $p>n, ~n\geq2$ and ${\mathcal{A}}>0$,
then there is no positive entire convex solution of \eqref{e1.1} for ideal domain $\Omega$.
\end{theorem}

In the definition of ideal domain, $\Omega_0$ is only bounded Lipschitz convex domain. Therefore, Lemma \ref{l3.1} can not be applied directly to $\Omega_0$. Fortunately, we have the following variant of Lemma \ref{l3.1}.

\begin{lemma}\label{l6.1}
Letting $\Omega$ be a bounded Lipschitz convex domain in ${\mathbb{R}}^n$ and $f(u)=u^p$ for some $p>n$, there exists a large convex solution $u\in C^\infty(\Omega)$
of \eqref{e1.1} satisfying \eqref{e3.2} for some positive constant $\mathcal{C}_3=\mathcal{C}_3(n,p,\Omega)$, where $\alpha=\frac{n+1}{p-n}$.
\end{lemma}

\begin{proof}
Assuming $\Omega$ is a bounded Lipschitz convex domain, there exists a sequence of smooth convex domains $\Omega_j\subset\Omega, j\in{\mathbb{N}}$ which is monotone increasing in $j$ and satisfies $\bigcup\limits_{j=1}^\infty\Omega_j=\Omega$. Applying Lemma \ref{l3.1} to $\Omega_j$, there exists a sequence of positive constants $C_{3,j}$ such that \eqref{e3.2} holds for solution $u_j$ of \eqref{e1.1} on $\Omega_j$. Noting that by the maximum principle of elliptic PDEs, we have $u_j$ is monotone decreasing in $j$ and $C_{3,j}$ can also be chosen to be monotone decreasing. As a result, passing to the limit $j\to\infty$ in \eqref{e3.2}, we obtain that for all $x\in\Omega$, the limiting solution $u(x):=\lim\limits_{j\to\infty}u_j(x)$ satisfying
\begin{equation*}
\begin{aligned}
u(x)&=\lim_{j\to\infty}u_j(x)\\
&\leq\lim_{j\to\infty}C_{3,j}\mathrm{dist}^{-\alpha}(x,\partial\Omega_j)\\
&\leq C_{3,\infty}\mathrm{dist}^{-\alpha}(x,\partial\Omega)
\end{aligned}
\end{equation*}
for some positive constant $C_{3,\infty}$. To show the left inequality of \eqref{e3.2}, one needs only to choose a monotone shrinking bounded smooth convex domains $\widetilde{\Omega}_j\supset\Omega$ satisfying $\bigcap\limits_{j=1}^\infty\widetilde{\Omega}_j=\Omega$. Then, applying the monotone increasing property of $\widetilde{u}_j$ and $\widetilde{C}_{3,j}^{-1}$, we arrive at
\begin{equation*}
\begin{aligned}
u(x)&\geq\lim_{j\to\infty}\widetilde{u}_j(x)\\
&\geq\lim_{j\to\infty}\widetilde{C}^{-1}_{3,j}\mathrm{dist}^{-\alpha}(x,\partial\Omega_j)\\
&\geq \widetilde{C}^{-1}_{3,\infty}\mathrm{dist}^{-\alpha}(x,\partial\Omega)
\end{aligned}
\end{equation*}
for any $x\in\Omega$ and some positive constant $\widetilde{C}_{3,\infty}$. Choosing
$$
\mathcal{C}_3:=\max\Big\{C_{3,\infty}, \widetilde{C}_{3,\infty}\Big\},
$$
we conclude that $u$ satisfies \eqref{e3.2} for $\mathcal{C}_3$.
\end{proof}

Now, we can complete the proof of Theorem \ref{t6.1}.

\medbreak
\noindent\textbf{Proof of theorem \ref{t6.1}}: As in the proof of Theorem \ref{t3.1}, let $u$ be a positive large convex solution of \eqref{e1.1} in $\Omega$.
Then, by assumption \eqref{e3.1}, $u$ is a supsolution to
\begin{equation}\label{e6.3}
\det D^2u={\mathcal{A}}u^p, ~~\forall x\in\Omega.
\end{equation}
Letting $\Omega_0\subset\Omega$, $a_k\in\Omega, \lambda_k\to\infty, k\in{\mathbb{N}}$ as given in definition of ideal domain, we set
$$
\Omega_k:=\lambda_k(\Omega-a_k)\subset\Omega, ~~\forall k\in{\mathbb{N}}.
$$
By Lemma \ref{l6.1}, we denote by $u_0$ a convex large solution of \eqref{e1.1} for $f(u)=u^p$ on $\Omega_0$, which satisfies that
$$
\mathcal{C}_3^{-1}\mathrm{dist}^{-\alpha}(x,\partial\Omega_0)\leq u_0(x)\leq \mathcal{C}_3\mathrm{dist}^{-\alpha}(x,\partial\Omega_0), ~~\forall x\in\Omega_0.
$$
Rescaling $u_0$ by
$$
u_0(x)=\lambda_k^{\frac{2n}{p-n}}u_k(\lambda_k(x-a_k)), ~~\forall x\in \Omega_0,
$$
then, for any $y\in\Omega_k$, the function $u_k(y)$ satisfies that
\begin{equation}\label{e6.4}
\mathcal{C}_3^{-1}\lambda_k^{\frac{n-1}{n-p}}\mathrm{dist}^{-\alpha}(y,\partial\Omega_k)\leq u_k(y)\leq
\mathcal{C}_3\lambda_k^{\frac{n-1}{n-p}}\mathrm{dist}^{-\alpha}(y,\partial\Omega_k)
\end{equation}
and
$$
\begin{cases}
\det D^2u_k=u_k^p, ~~&\text{in} ~~\Omega_k,\\
u_k(y)\to\infty, ~~&\text{as}~~y\to\partial\Omega_k.
\end{cases}
$$
Letting $w_k(y):={\mathcal{A}}^{\frac{1}{n-p}}u_k(y)$, then $w_k$ is a solution to
$$
\begin{cases}
\det D^2w_k={\mathcal{A}}w_k^p, ~~&\text{in} ~~\Omega_k,\\
w_k(y)\to\infty, ~~&\text{as}~~y\to\partial\Omega_k.
\end{cases}
$$
Next, let us denote by $\left[U^{ij}(t)\right]$ the cofactor matrix of $tD^2w_{k}+(1-t)D^2u$ with $t\in [0,1]$. By some simple calculations, it is no difficult to find that
$$
\begin{cases}
\displaystyle \int^1_0U^{ij}(t)dt\cdot D_{ij}(w_{k}-u)\le {\mathcal{B}}(w_{k}-u), & \text{in}~~\Omega_k,\\
w_{k}-u\geq0, & \text{on}~~\partial \Omega_k,
\end{cases}
$$
where the coefficient $\mathcal{B}$ is defined as
$$
{\mathcal{B}}=\frac{f(w_{k})-f(u)}{w_{k}-u}\geq 0.
$$
It is inferred from the maximum principle of elliptic PDEs that $u\leq w_{k}$ in $\Omega_k$. Combining \eqref{e6.4} and the definition of $w_{k}$, we immediately obtain
\begin{equation}\label{e6.5}
u(y)\leq \mathcal{C}_3{\mathcal{A}}^{\frac{1}{n-p}} \lambda_k^{\frac{n-1}{n-p}}\mathrm{dist}^{-\alpha}(y,\partial \Omega_k), \ \ \forall y\in \Omega_k.
\end{equation}
Taking $y$ to be an ideal center $a_\infty$ of $\Omega$ and sending $k$ to be large, it is inferred from \eqref{e6.5} that $u(a_\infty)$ must be identical to zero since $p>n$, which contradicts with the positivity of $u$. The proof of the theorem is complete.\hfill\ensuremath{\square}

\bigbreak
\section{Wiping barrier and large solution on unbounded domain}

\bigbreak
\noindent Inspired by Theorems \ref{t3.1}, \ref{t4.1} and \ref{t6.1}, it is natural to ask the following question.

\medbreak
\noindent\textbf{Question:} For $n\geq2$ and $p<n$, given any unbounded convex domain $\Omega\not={\mathbb{R}}^n$, whether there is some Euclidean complete solution of \eqref{e0.1} on $\Omega$?

\medbreak
It would be surprising that we have the following non-existence result of large solution of \eqref{e0.1}.

\begin{theorem}\label{t7.1}
Let $n=2$, $p\in(0,1/2)$, and $\Omega\neq\mathbb{R}^2$ be an unbounded domain. Then there is no positive convex Euclidean complete solution to \eqref{e0.1}.
\end{theorem}

Our proof of Theorem \ref{t7.1} relies on a construction of zero barrier solution on some wiping domain defined by
$$
\Sigma_{\beta,r_0}:=\Big\{(x,y)\in{\mathbb{R}}^2~\Big|~e^xy^\beta\in(r_0,\infty)\Big\}
$$
where $r_0\in(0,\infty)$ and $\beta<0$. We need to use the following existence result in the proof of Theorem \ref{t7.1}.

\begin{proposition}\label{p7.1}
Considering \eqref{e0.1} for $n=2, ~p\in(0,1/2)$ and $\beta<0$, there exists a positive constant $r_0=r_0(p,\beta)$ such that \eqref{e0.1} admits a positive smooth convex solution $u$ on $\Sigma_{\beta,r_0}$ satisfying $u(x,y)=0$ for any $(x,y)\in \varsigma_{\beta,r_0}$, where
$$
\varsigma_{\beta,r_0}:=\Big\{(x,y)\in{\mathbb{R}}^2~\Big|~e^xy^\beta=r_0\Big\}.
$$
\end{proposition}

To verify the validity of Proposition \ref{p7.1}, we look for solutions of the form
$$
u(x,y)=y^\alpha\varphi(e^xy^\beta),~~x\in{\mathbb{R}},~y\in{\mathbb{R}}_+.
$$
Letting $r:=e^xy^{\beta}$, direct computation shows that
\begin{equation*}
\begin{aligned}
u_{xx}&=y^{\alpha}\Big[r\varphi_{r}+r^2\varphi_{rr}\Big],\\
u_{xy}&=y^{\alpha-1}\Big[(\alpha+\beta)r\varphi_r+\beta r^2\varphi_{rr}\Big],\\
u_{yy}&=y^{\alpha-2}\Big[\alpha(\alpha-1)\varphi+\beta(2\alpha+\beta-1)r\varphi_r+\beta^2r^2\varphi_{rr}\Big].
\end{aligned}
\end{equation*}
Therefore, we can obtain that
\begin{equation*}
\begin{aligned}
\det D^2u&=y^{2\alpha-2}\Big[r^2\varphi_{rr}\mathcal{P}_1(r)+\mathcal{P}_2(r)\Big],
\end{aligned}
\end{equation*}
where
\begin{equation*}
\begin{aligned}
&\mathcal{P}_1(r):=\alpha(\alpha-1)\varphi-\beta r\varphi_r,\\
&\mathcal{P}_2(r):=\alpha(\alpha-1)r\varphi\varphi_r-(\alpha^2+\beta)r^2\varphi_r^2.
\end{aligned}
\end{equation*}
Noticing that $u^p=y^{\alpha p}\varphi^p$, if one chooses $\alpha=\frac{2}{2-p}$, then \eqref{e0.1} reduces to
\begin{equation}\label{e7.1}
\begin{aligned}
r^2\varphi_{rr}\mathcal{P}_1(r)+\mathcal{P}_2(r)=\varphi^p.
\end{aligned}
\end{equation}
Setting $\xi(\varphi)=r\varphi_r$, one has
$$
\xi\xi'=\xi\frac{d\xi}{d\varphi}=r\varphi_r+r^2\varphi_{rr}=\xi+r^2\varphi_{rr},
$$
and thus we can rewrite \eqref{e7.1} as
\begin{equation}\label{e7.2}
\xi\xi'=\frac{\alpha^2\xi^2+\varphi^p}{\alpha(\alpha-1)\varphi-\beta\xi}.
\end{equation}
To recover $\varphi$ from $\xi$, we choose $\varphi_1=\varphi(r_1)$ for some $r_1>0$, and obtain that
\begin{equation}\label{e7.3}
\int^{\varphi(r)}_{\varphi_1}\frac{d\varphi}{\xi(\varphi)}=\int^r_{r_1}\frac{dr}{r}=\log\frac{r}{r_1}.
\end{equation}
Now, let us introduce a closed convex subset
\begin{eqnarray*}
\mathcal{F}_{\beta,q,\delta}:=\Big\{\xi\in C([0,\delta])~\Big|~ \xi(\varphi)-\gamma_\beta\varphi^{\frac{p+1}{3}}\in[-\varphi^{q},\varphi^{q}], \ \ \forall\varphi\in[0,\delta]\Big\}
\end{eqnarray*}
of $C([0,\delta])$, where $q\in(\frac{p+1}{3},1)$ is chosen later and
\begin{equation}\label{e7.4}
\gamma_\beta:=\left[\frac{3}{|\beta|(p+1)}\right]^{\frac{1}{3}}.
\end{equation}
For each $\xi\in\mathcal{F}_{\beta,q,\delta}$, we also define a mapping $T\xi=:\zeta$ by
\begin{equation}\label{e7.5}
\begin{cases}
\displaystyle \zeta'=\frac{\alpha^2\xi+\varphi^p\xi^{-1}}{\alpha(\alpha-1)\varphi-\beta\xi},~~\varphi\in (0,\delta], \\
\zeta(0)=0.
\end{cases}
\end{equation}

\begin{lemma}\label{l7.1}
Supposing that $p\in(0,1/2)$ and $\beta<0$, if one chooses $q$ closing to 1 then chooses $\delta$ small, then $T$ is a continuous and compact mapping from $\mathcal{F}_{\beta,q,\delta}$ to $\mathcal{F}_{\beta,q,\delta}$.
\end{lemma}

\begin{proof}
To verify the compactness of $T$, we use \eqref{e7.5} to calculate
\begin{equation*}
\begin{aligned}
\zeta'&\leq\frac{\alpha^2\Big(\gamma_\beta\varphi^{\frac{p+1}{3}}+\varphi^q\Big)+\varphi^p\Big(\gamma_\beta\varphi^{\frac{p+1}{3}}-\varphi^q\Big)^{-1}}{\alpha(\alpha-1)\varphi-\beta\Big(\gamma_\beta\varphi^{\frac{p+1}{3}}-\varphi^q\Big)}\\
&\leq|\beta|^{-1}\gamma_\beta^{-2}\varphi^{\frac{p-2}{3}}+2|\beta|^{-1}\gamma_\beta^{-3}\varphi^{q-1}+o(\varphi^{q-1})\\
&\leq\frac{\gamma_\beta(p+1)}{3}\varphi^{\frac{p-2}{3}}+q\varphi^{q-1}, \ \ \forall\varphi\in(0,\delta]
\end{aligned}
\end{equation*}
and
\begin{equation*}
\begin{aligned} \zeta'&\geq\frac{\alpha^2\Big(\gamma_\beta\varphi^{\frac{p+1}{3}}-\varphi^q\Big)+\varphi^p\Big(\gamma_\beta\varphi^{\frac{p+1}{3}}+\varphi^q\Big)^{-1}}{\alpha(\alpha-1)\varphi-\beta\Big(\gamma_\beta\varphi^{\frac{p+1}{3}}+\varphi^q\Big)}\\
&\geq|\beta|^{-1}\gamma_\beta^{-2}\varphi^{\frac{p-2}{3}}-2|\beta|^{-1}\gamma_\beta^{-3}\varphi^{q-1}+o(\varphi^{q-1})\\
&\geq\frac{\gamma_\beta(p+1)}{3}\varphi^{\frac{p-2}{3}}-q\varphi^{q-1}, \ \ \forall\varphi\in(0,\delta]
\end{aligned}
\end{equation*}
provided $q$ is chosen to close 1 and then $\delta$ is chosen small, where
$$
2|\beta|^{-1}\gamma_\beta^{-3}=\frac{2(p+1)}{3}<1
$$
has been used. Integrating over $\varphi$ yields that
$$
\zeta(\varphi)-\gamma_\beta\varphi^{\frac{p+1}{3}}\in[-\varphi^q,\varphi^q], \ \ \forall\varphi\in(0,\delta].
$$
Consequently, $T$ is a continuous and compact mapping from $\mathcal{F}_{\beta,q,\delta}$ to $\mathcal{F}_{\beta,q,\delta}$. This completes the proof.
\end{proof}

As a corollary, we obtain the following local solvability of \eqref{e7.2}.

\begin{corollary}\label{c7.1}
Under the assumptions of Lemma \ref{l7.1}, the mapping $T$ admits a fixed point $\zeta$ in $\mathcal{F}_{\beta,q,\delta}$ which is a local smooth solution of \eqref{e7.2}.
\end{corollary}

Finally, let us complete the proof of Proposition \ref{p7.1} with the help of the following long time existence and asymptotic results.

\begin{proposition}\label{p7.2} {\rm(Long time existence)}
Under the assumptions of Corollary \ref{c7.1}, the local solution $\zeta$ exists for all time $\varphi>0$ and preserves positivity and monotonicity.
\end{proposition}

\begin{proof}
It is clear that the solution of \eqref{e7.2} preserves positivity and monotonicity. To show that the solution exists for all $\varphi$, one needs only using\eqref{e7.2} to deduce that
\begin{equation*}
\begin{aligned}
(\zeta^2)'\leq\frac{2\alpha}{\alpha-1}\varphi^{-1}\zeta^2+\frac{2}{\alpha(\alpha-1)}\varphi^{p-1}.
\end{aligned}
\end{equation*}
This is equivalent to
\begin{equation*}
\Big(\varphi^{-\frac{2\alpha}{\alpha-1}}\zeta^2\Big)'\leq\frac{2}{\alpha(\alpha-1)}\varphi^{p-1-\frac{2\alpha}{\alpha-1}}.
\end{equation*}
Integrating the above inequality immediately yields
\begin{equation*}
\varphi^{-\frac{2\alpha}{\alpha-1}}\zeta^2\leq\delta^{-\frac{2\alpha}{\alpha-1}}                                       \zeta^2(\delta)-\frac{2p}{\alpha(\alpha-1)(4-p^2)}\Big(\varphi^{p-\frac{2\alpha}{\alpha-1}}-\delta^{p-\frac{2\alpha}{\alpha-1}}\Big)
\end{equation*}
for all $\varphi\ge \delta$. Hence, $\zeta$ is A-priori bounded from above and thus exists for all $\varphi$.
\end{proof}

The second proposition gives the asymptotic behavior of $\zeta$ at infinity.

\begin{proposition}\label{p7.3} {\rm(Asymptotic behavior)}
Under the assumptions of Corollary \ref{c7.1}, for each $\varepsilon>0$, there exists a positive constant $\varphi_\varepsilon$ such that
\begin{equation}\label{e7.6}
\left(\frac{\alpha}{|\beta|}-\varepsilon\right)\varphi\leq \zeta\leq\left(\frac{\alpha}{|\beta|}+\varepsilon\right)\varphi,~~
\forall\varphi\geq\varphi_\varepsilon
\end{equation}
holds.
\end{proposition}

We divided the proof into two lemmas.

\begin{lemma}\label{l7.2}
Under the assumptions of Corollary \ref{c7.1}, there holds
\begin{equation}\label{e7.7}
\zeta^2\leq A_1\varphi^{\frac{4}{p}}-A_2\varphi^p, ~~\forall\varphi\geq\delta,
\end{equation}
where
\begin{equation*}
\begin{aligned}
&A_1:=\delta^{-\frac{4}{p}}\zeta^2(\delta)+\frac{2p}{\alpha(\alpha-1)(4-p^2)}\delta^{p-\frac{4}{p}}\\
&A_2:=\frac{2p}{\alpha(\alpha-1)(4-p^2)}.
\end{aligned}
\end{equation*}
\end{lemma}

Lemma \ref{l7.2} is a direct consequence of the proof of Proposition \ref{p7.2}.

\begin{lemma}\label{l7.3}
Under the assumptions of Corollary \ref{c7.1}, for any $\varepsilon>0$, there exists positive constant $\varphi_2=\varphi_2(\varepsilon)$ such that
\begin{equation}\label{e7.8}
\zeta\geq\left(\frac{\alpha^2(2-p)}{2|\beta|}-\varepsilon\right)\varphi, ~~\forall\varphi\geq\varphi_2.
\end{equation}
\end{lemma}

\begin{proof}
By \eqref{e7.2} and $p<2$, there holds
\begin{equation*}
\begin{aligned}
\zeta'\geq\frac{\alpha^2\zeta}{\alpha(\alpha-1)\varphi-\beta\zeta}
&\Rightarrow\frac{d\varphi}{d\zeta}\leq\frac{\alpha-1}{\alpha}\zeta^{-1}\varphi+\frac{|\beta|}{\alpha^2}\\
&\Rightarrow\Big(\zeta^{-\frac{\alpha-1}{\alpha}}\varphi\Big)'\leq\frac{|\beta|}{\alpha^2}\zeta^{-\frac{\alpha-1}{\alpha}}\\
&\Rightarrow\zeta^{-\frac{p}{2}}\varphi\leq C_{11}+\frac{2|\beta|}{\alpha^2(2-p)}\zeta^{1-\frac{p}{2}}
\end{aligned}
\end{equation*} for some universal constant $C_{11}>0$. Then, by Young's inequality, it follows that
\begin{equation*}
\begin{aligned}
\varphi&\leq C_{11}\zeta^{\frac{p}{2}}+\frac{2|\beta|}{\alpha^2(2-p)}\zeta\\
&\leq\left(\frac{2|\beta|}{\alpha^2(2-p)}+\varepsilon\right)\zeta+C_\varepsilon
\end{aligned}
\end{equation*}
for each $\varepsilon>0$, where $C_\varepsilon$ is a positive constant depending on $\varepsilon$. Therefore, for this given $\varepsilon$, there exists $\varphi_2=\varphi_2(\varepsilon)$ such that \eqref{e7.8} holds.
\end{proof}

Now, one can complete the proof of Proposition \ref{p7.3}.

\medbreak
\noindent\textbf{Proof of Proposition \ref{p7.3}.} By \eqref{e7.2} and Lemma \ref{l7.3}, for each $\varepsilon$ small, there exists $\varphi_3=\varphi_3(\varepsilon)$ such that
\begin{equation}\label{e7.9}
\zeta'\leq\frac{(\alpha^2+\varepsilon)\zeta}{\alpha(\alpha-1)\varphi-\beta\zeta}, ~~\forall\varphi\geq\varphi_3.
\end{equation}
Solving this first order O.D.E. as in the proof of Lemma \ref{l7.3}, it yields the desired inequality \eqref{e7.6} with the help of \eqref{e7.8}.\hfill\ensuremath{\square}

\medbreak
At the end of this section, let us complete the proof of Proposition \ref{p7.1} and Theorem \ref{t7.1}.

\medbreak
\noindent\textbf{Proof of Proposition \ref{p7.1}.} Noting that the solution $\zeta$ derived by Corollary \ref{c7.1} and Proposition \ref{p7.2} satisfies
$$
\zeta(\varphi)\sim\gamma_\beta\varphi^{\frac{p+1}{3}}, ~~\forall\varphi\sim0,
$$
there exists a positive constant $r_0\in(0,r_1)$, such that
\begin{equation}\label{e7.10}
\int^0_{\varphi_1}\frac{d\varphi}{\zeta(\varphi)}=\log\frac{r_0}{r_1}.
\end{equation}
Similarly, by Proposition \ref{p7.3},
$$
\zeta(\varphi)\sim\left(\frac{\alpha}{|\beta|}\pm\varepsilon\right)\varphi, ~~\forall\varphi\sim\infty,
$$
there holds
\begin{equation}\label{e7.11}
\int^\infty_{\varphi_1}\frac{d\varphi}{\zeta(\varphi)}=\infty.
\end{equation}
Therefore, the recovered convex solution $\varphi$ from \eqref{e7.3} satisfies that
\begin{equation}\label{e7.12}
\varphi(r_0)=0, ~~\varphi(\infty)=\infty.
\end{equation}
Using the relation $u(x,y)=y^\alpha\varphi(e^xy^\beta)$ and $\alpha>0, \beta<0$, one obtains a desired positive convex smooth solution $u$ of $\Sigma_{\beta,r_0}$. The proof of Proposition \ref{p7.1} is complete.\hfill\ensuremath{\square}

\medbreak
\noindent\textbf{Proof of Theorem \ref{t7.1}.} Suppose on the contrary, there exists a positive convex large solution $u$ of \eqref{e0.1} on an unbounded convex domain $\Omega\not={\mathbb{R}}^2$. Choosing a point $z_0=(x_0,y_0)\in\partial\Omega$ and drawing a line $l_{z_0}$ tangential to $\partial\Omega$ at the point $z_0$, whose unit normal is given by $\nu$. The boundary of $\Sigma_{\beta,r_0}$ contains a curve $\varsigma_{\beta,r_0}$ and a straight line $l_{\beta,r_0}$. The curve $\varsigma_{\beta,r_0}$ divides the plane into two sides. The side of $\varsigma_{\beta,r_0}$ which contains $\Sigma_{\beta,r_0}$ will be called positive side of $\varsigma_{\beta,r_0}$ for short. Now, translating and rotating the wiping domain $\Sigma_{\beta,r_0}$ obtained in Proposition \ref{p7.1}, $\Omega$ is divided into two portions by $\varsigma_{\beta,r_0}$. One may assume that one of the portions $\Sigma_*$, which lies on positive side of $\varsigma_{\beta,r_0}$, is bounded and lies strictly inside of $\Sigma_{\beta,r_0}$. Comparing with the zero barrier solution $u_*$ found in Proposition \ref{p7.1} on $\Sigma_*$, one gets
\begin{equation}\label{e7.13}
u(x,y)\geq u_*(x,y)/\varepsilon, ~~\forall(x,y)\in\Sigma_*
\end{equation}
for each $\varepsilon\in(0,1)$, where we have used $u_*/\varepsilon$ is a subsolution of \eqref{e0.1} and
\begin{equation}\label{e7.14}
u(x,y)\geq u_*(x,y)/\varepsilon, ~~\forall(x,y)\in\partial\Sigma_*.
\end{equation}
Finally, let $\varepsilon\to 0^+$ in \eqref{e7.13}, then $u$ must be infinite everywhere in $\Sigma_*$. Contradiction holds. The proof of Theorem \ref{t7.1} is complete. \hfill\ensuremath{\square}

\bigbreak
{\noindent \textbf{Acknowledgements:}}

\medbreak
The author (SZ) would like to express his deepest gratitude to Professors Xi-Ping Zhu, Kai-Seng Chou, Xu-Jia Wang and Neil Trudinger for their constant encouragements and warm-hearted helps. This paper is also dedicated to the memory of Professor Dong-Gao Deng.


\vspace{15pt}
\leftline{\bfseries \large References}
\begin{thebibliography}{1} \small

\bibitem{BM1} C. Bandle and M. Marcus, {\it Large solutions of semilinear elliptic equations: existence, uniqueness and asymptotic behavior}, J. Anal. Math. {\bf58} (1992), 9-24.

\bibitem{BM2} C. Bandle and M. Marcus, {\it Asymptotic behaviour of solutions and their derivatives, for semilinear elliptic problems with blowup on the boundary}, Ann. Inst. H. Poincar\'{e} Anal. Non Lin\'{e}aire {\bf12} (1996), 155-171.

\bibitem{BCGJ} J.G. Bao, J.Y. Chen, B. Guan and M. Ji, {\it Liouville property and regularity of a Hessian quotient equation}, Amer. J. Math. {\bf 125} (2003), 301-316.

\bibitem{Bie1} L. Bieberbach, {\it $\Delta u=e^u$ und die automorphen Funktionen}, Math. Ann., {\bf77} (1916), 173-212.

\bibitem{CSF} A. Colesanti, P. Salani and E. Francini, {\it Convexity and asymptotic estimates for large solutions of Hessian equations}, Differ. Integral Equ. {\bf13} (2000), 1459-1472

\bibitem{CD} F.C. C\^{i}rstea and Y. Du, {\it General uniqueness results and variation speed for blow-up solutions of elliptic equations}. Proc. Lond. Math. Soc., {\bf91} (2005), 459-482.

\bibitem{CR1} F.C. C\^{i}rstea and V. R\u{a}dulescu, {\it Uniqueness of the blow-up boundary solution of logistic equations with absorption}, C. R. Math. Acad. Sci. Paris {\bf335} (2002), 447-452.

\bibitem{CR2} F.C. C\^{i}rstea and V. R\u{a}dulescu, {\it Nonlinear problems with boundary blow-up: a Karamata regular variation theory approach}, Asymptot. Anal. {\bf46} (2006), 275-298.

\bibitem{CT} F.C. Cirstea and C. Trombetti, {\it On the Monge-Amp\`{e}re equation with boundary blow-up: existence, uniqueness and asymptotics}, Calc. Var., {\bf31}
(2008), 167-186.

\bibitem{CN1} K.S. Cheng and W.M. Ni, {\it On the structure of the comformal scalar curvature equation on $\mathbb{R}^n$}, Indiana Univ. Math. J., {\bf41} (1992), 261-278.

\bibitem{CW1} K.S. Chou and X.J. Wang, {\it A variational theory of the Hessian equation}, Comm. Pure Appl. Math., {\bf54} (2002), 1029-1064.

\bibitem{CW2} K.S. Chou and X.J. Wang, {\it Entire solutions of the Monge-Amp\`{e}re equation}, Comm. Pure Appl. Math., {\bf49} (1996), 529-539.

\bibitem{Caf1} L.A. Caffarelli, {A localization property of viscosity solutions to the Monge-Amp\`{e}re equation and their
strict convexity}, Ann. Math. {\bf131} (1990), 129-134.

\bibitem{CNS} L.A. Caffarelli, L. Nirenberg and J. Spruck, {\it The Dirichlet problem for nonlinear second-order elliptic equa
tions. I. Monge-Amp\`{e}re equation}, Comm. Pure Appl. Math. {\bf37} (1984), 369-402.

\bibitem{CY1} S.Y. Cheng and S.T. Yau, {\it On the existence of a complete K\"{a}hler metric on noncompact complex manifolds and the regularity of Fefferman's
equation}, Comm. Pure Appl. Math., {\bf33} (1980), 507-554.

\bibitem{CY2} S.Y. Cheng and S.T. Yau, {\it The real Monge-Amp\`{e}re equation and affine flat structures. In: Chern, S.S., Wu, W. (eds.)} Proceedings of 1980 Beijing
Symposium on Differential Geometry and Differential Equations, vol.1, 339-370, Beijing. Science Press, New York (1982).

\bibitem{Du} S.Z. Du, {\it Bernstein problem of affine maximal type hypersurfaces on dimension $N\ge 3$}, J. Differential Equations, {\bf269} (2020), 7429-7469.

\bibitem{FZ} M.Q. Feng and X.M. Zhang, {\it Strictly convex solutions to the singular boundary blow-up Monge-Amp\`{e}re problems: existence and asymptotic behavior}, J. Gemo. Anal., {\bf 34} (2024), 301.

\bibitem{Guan} B. Guan, {\it The Dirichlet problem for a class of fully nonlinear elliptic equations}, Comm. Partial Differ. Equ., {\bf19} (1994), 399-416.

\bibitem{GJ} B. Guan and H. Jian, {\it The Monge-Amp\`{e}re equation with infinite boundary value}, Pac. J. Math., {\bf216} (2004), 77-94.


\bibitem{Huang1} Y. Huang, {\it Boundary asymptotical behavior of large solutions to Hessian equations}, Pacific J. Math., {\bf244} (2010), 85-98.

\bibitem{JW1} H.Y. Jian and X.J. Wang, {\it Existence of entire solutions to the Monge-Amp\`{e}re equation}, Amer. J. Math., {\bf136} (2014), 1093–1106.

\bibitem{JLX1} Q.N. Jin, Y.Y. Li and H.Y. Xu, {Nonexistence of positive solutions for some fully nonlinear elliptic equations}, Methods Appl. Anal. {\bf12} (2005),  441-449.

\bibitem{JB1} X.H. Ji and J.G. Bao, {\it Necessary and sufficient conditions on solvability for Hessian inequalities}, Proc. Amer. Math. Soc., {\bf138} (2010), 175-188.

\bibitem{Kel1} J.B. Keller, {\it On solutions of $\Delta u=f(u)$}, Comm. Pure Appl. Math., {\bf10} (1957), 503-510.

\bibitem{KN1} V.A. Kondrat\'{e}v and V.A. Nikishkin, {\it Asymptotics near the boundary of a solution of a singular boundary-value problem for a semilinear elliptic equation}, Differ. Equ., {\bf26} (1990), 345-348.

\bibitem{LM1} A.C. Lazer and P.J. McKenna, {\it On singular boundary value problems for the Monge-Amp\`{e}re operator}, J. Math. Anal. Appl., {\bf197} (1996),
341-362.

\bibitem{LN1} C. Loewner and L. Nirenberg, {\it Partial differential equations invariant under conformal or projective transformations, in:
Contributions to Analysis (A Collection of Papers Dedicated to Lipman Bers)}, Academic Press, New York, 1974, pp. 245-272.

\bibitem{Lion1} P.L. Lions, {\it Sur les \'{e}quations de Monge-Amp\`{e}re equations. (French) [On Monge-Amp\`{e}re equations]}, Arch. Ration. Mech. Anal., {\bf 89} (1985), 93-122.

\bibitem{Lion2} P.L. Lions, {\it Two remarks on Monge-Amp\`{e}re equations}, Ann. Mat. Pura Appl., {\bf142} (1985), 263-275.

\bibitem{LB1} X. Li and J.G. Bao, {\it Existence and asymptotic behavior of entire large solutions for Hessian equations}, Commun. Pure Appl. Anal., {\bf23} (2024), 253-268.

\bibitem{Moh1} A. Mohammed, {\it On the existence of solutions ot the Monge-Amp\`{e}re equation with infinite boundary values}, Proc. Amer. Math. Soc., {\bf135} (2007),
141-149.

\bibitem{MV1} M. Marcus and L. V\'{e}ron, {\it Uniqueness and asymptotic behaviour of solutions with boundary blow-up for a class of nonlinear elliptic equations}. Ann. Inst. H. Poincar\'{e} Anal. Non Lin\'{e}aire, {\bf14} (1997), 237-275.

\bibitem{MV2} M. Marcus and L. V\'{e}ron, {\it Existence and uniqueness results for large solutions of general nonlinear elliptic equations}. J. Evol. Equ.,
{\bf3} (2003), 637-652.

\bibitem{Mat1} J. Matero, {\it The Bieberbach-Rademacher problem for the Monge-Amp\`{e}re operator}, Manuscr. Math., {\bf91} (1996), 379-391.

\bibitem{ML} S.S. Ma and D.S. Li, {\it Exstence and boundary asymptotic behavior of large solutions of Hessian equations}, Nonlinear Analysis, {\bf187} (2019),
1-17.

\bibitem{Oss1} R. Osserman, {\it On the inequality $\Delta u \ge f(u)$}, Pacific J. Math., {\bf7} (1957), 1641-1647.

\bibitem{Pog} A.V. Pogorelov, {\it The multidimensional Minkovski problem}, Wiley, New York, (1978).

\bibitem{Pli} S. Pli\'{s}, {\it On boundary blow-up problems for the complex Monge-Amp\`{e}re equation}, Proc. Am. Math. Soc., {\bf136} (2008), 4355-4364.

\bibitem{Rad1} H. Rademacher, {\it Einige besondere probleme partieller Differentialgleichungen. In: Die Differential und Integralgleichungen der Mechanik und Physik}, I, Rosenberg, New York, 2nd edn, (1943), pp. 838-845

\bibitem{Sal1} P. Salani, {\it Boundary blow-up problems for Hessian equations}, Manuscripta Math., {\bf96} (1998), 281-294.

\bibitem{Tak1} K. Takimoto, {\it Solution to the boundary blow-up problem for $k$-curvature equation}, Calc. Var. Partial Differ. Equ., {\bf26} (2006), 357-377.

\bibitem{Tso} K. Tso, {\it On a real Monge-Amp\`{e}re functional}, Invent. Math., {\bf101} (1990), 425-448.

\bibitem{Tru1} N.S. Trudinger, {\it Fully nonlinear, uniformly elliptic equations under natural structure conditions}, Trans. Am. Math. Soc., {\bf278} (1983), 751-769.

\bibitem{Tru2} N.S. Trudinger, {\it On the Dirichlet problem for Hessian equations}, Acta Math., {\bf175} (1995), 151-164.

\bibitem{Tru3}  N.S. Trudinger, {\it Weak solutions of Hessian equations}, Comm. Partial Differ. Equ. {\bf22} (1997) 1251-1261.

\bibitem{TU} N.S. Trudinger and J. Urbas, {\it The Dirichlet problem for the equation of prescribed Gauss curvature}, Bull. Aust. Math. Soc., {\bf28} (1983), 217-231.

\bibitem{TW} N.S. Trudinger and X.J Wang, {\it Hessian measure \MakeUppercase{\romannumeral 2}}, Ann. Math., {\bf 19} (1999), 579-604.

\bibitem{Wang1} X.J. Wang, {\it A class of fully nonlinear elliptic equations and related functionals}, Indiana Univ. Math. J, {\bf43} (1994), 25-54.

\bibitem{YC} H.T. Yang and Y.B. Chang, {On the blow-up boundary solutions of the Monge-Amp\`{e}re equation with singular weights}, Commun. Pure Appl. Anal., {\bf11} (2012), 697-708.

\bibitem{ZD} X.M. Zhang and Y.H. Du, {\it Sharp conditions for the existence of boundary blow-up solutions to the Monge-Amp\`{e}re equation}, Calc. Var. Partial Differ. Equ., {\bf57} (2018), 1-24.

\bibitem{ZF} X.M. Zhang and M.Q. Feng, {\it The existence and asymptotic behavior of boundary blow-up solutions to the $k$-Hessian equation}, J. Differential Equations, {\bf267} (2019), 4626-4672.

\bibitem{ZF2} X.M. Zhang and M.Q. Feng, {\it Blow-up solutions to the Monge–Amp\`{e}re equation with a gradient term: sharp conditions for the existence and asymptotic estimates}, Calc. Var. Partial Differ. Equ., {\bf61}(2022), 208.

\bibitem{Zhang1} Z.J. Zhang, {\it Boundary behavior of large solutions to the Monge-Amp\`{e}re equations with weights}, J. Differential Equations, {\bf259} (2015),
2080-2100.

\bibitem{Zhang2} Z.J. Zhang, {\it Large solutions to the Monge-Amp\`{e}re equations with nonlinear gradient terms: existence and boundary behavior}, J. Differential Equations, {\bf264} (2018), 263-296.

\bibitem{Zhang3} Z.J. Zhang, {\it Boundary behavior of large solutions to the Monge-Amp\`{e}re equation in a borderline case}, Acta Math. Sin. Engl. Ser., {\bf35} (2019), 1190-1204.

\bibitem{Zhang4} Z.J. Zhang, {\it Optimal global and boundary asymptotic behavior of large solutions to the Monge-Amp\`{e}re equation}, J. Funct. Anal., {\bf278} (2020), 108512.

\bibitem{Zhang5} Z.J. Zhang, {\it Boundary behavior of large solutions for equations of Monge-Amp\`{e}re type}, Chin. Ann. Math. Ser. B, (2025).

\bibitem{ZMML} Z.J. Zhang, Y.J. Ma, L. Mi and X.H. Li, {\it Blow-up rates of large solutions for elliptic equations}, J. Differential Equations, {\bf249} (2010), 180-199.

\bibitem{ZZ1} Z.J. Zhang and S. Zhou, {\it Existence of entire positive $k$-convex radial solutions to Hessian equations and systems with weights}., Appl. Math. Lett. {\bf50} (2015), 48-55.
\end{thebibliography}
\end{document}